\documentclass[english,11pt,a4paper,twoside,reqno]{article}

\usepackage{amsfonts}
\usepackage{amsmath}
\usepackage{amssymb}
\usepackage{textcomp}
\usepackage{mathrsfs}
\usepackage[T1]{fontenc}
\usepackage{graphics}
\usepackage{color}
\usepackage[english]{babel}
\usepackage[active]{srcltx}
\usepackage{bm}
\usepackage{enumerate}
\usepackage{eufrak}
\usepackage{euscript}
\numberwithin{equation}{section}
\usepackage{newlfont}
\usepackage{latexsym}
\usepackage{graphicx}
\usepackage{a4wide}

\usepackage{cite}        
\usepackage{url,float}         

\title{A uniform Tauberian theorem in optimal control}
\author{Miquel Oliu Barton\footnote{\'Equipe Combinatoire et Optimisation, CNRS FRE3232, Universit\'e Paris 6, UFR 929, 
 Paris, France.}
\ and Guillaume Vigeral\footnote{INRIA Saclay \& Centre de
Math\'ematiques Appliqu\'ees, \'Ecole Polytechnique, 91128 Palaiseau,
France.}}

\date{\today}

\newcounter{compteur}

\newcounter{propo}

\newtheorem{prop}[propo]{Proposition}
\newtheorem{theo}[compteur]{Theorem}
\newtheorem{lemme}[propo]{Lemma}

\newtheorem{remarque}{Remark}

\newenvironment{proof}[1][Proof]{\noindent \textbf{#1.}~ }
{\hfill\rule{2mm}{2mm} \vspace{\parskip} }

\begin{document}
 \maketitle

\textbf{Abstract} : In an optimal control framework, we consider the
value $V_T(x)$ of the problem starting from state $x$ with finite
horizon $T$, as well as the value $V_\lambda(x)$ of the
$\lambda$-discounted problem starting from $x$. We prove that
uniform convergence (on the set of states) of the values
$V_T(\cdot)$ as $T$ tends to infinity is equivalent to uniform
convergence of the values $V_\lambda(\cdot)$ as $\lambda$ tends to
0, and that the limits are identical. An example is also provided to
show that the result does not hold for pointwise convergence. This
work is an extension, using similar techniques, of a related result in
a discrete-time framework \cite{LehSys}.

\section{Introduction}\label{Controleoptimal}
Finite horizon problem of optimal control have been studied
intensively since the pioneer work of
Stekhov, Pontryagin, Boltyanski, Hestenes, Bellman and Isaacs during
the cold war - see for instance \cite{LeeMarkus,Kirk,Bardi} for major references,
or \cite{Evans} for a short, clear introduction. A classical model considers the following controlled dynamic
over $\mathbb{R}+$
\begin{equation}\label{eqcontrole} \begin{cases}
y'(s)=f(y(s),u(s))\\
y(0)=y_0
\end{cases}
\end{equation}
\noindent where $y$ is a function from $\mathbb{R}+$ to
$\mathbb{R}^n$, $y_0$ is a point in $\mathbb{R}^n$, $u$ is the
control function which belongs to $\mathcal{U}$, the set of
Lebesgue-measurable functions from $\mathbb{R}+$ to a metric space
$U$ and the function $f:\mathbb{R}^n \times U \to \mathbb{R}^n$
satisfies the usual conditions, that is: Lipschitz with respect to
the state variable, continuous with respect to the control variable
and bounded by a linear function of the state variable, for any
control $u$.

Together with the dynamic, an objective function $g$ is given, interpreted as the cost function which is to be minimized and assumed to be Borel-measurable
from $\mathbb{R}^n \times U$ to $[0,1]$. For each finite horizon
$t\in ]0,+\infty[$, the average value of the optimal control problem with horizon $t$ is defined as
\begin{equation}
V_t(y_0)=\inf_{u\in\mathcal{U}}\frac{1}{t}\int_0^t
g(y(s,u,y_0),u(s))ds \end{equation}
 It is quite natural to define,
whenever the trajectories considered are infinite, for any discount
factor $\lambda >0$, the $\lambda$-discounted value of the optimal
control problem, as
\begin{equation} V_\lambda(y_0)=\inf_{u\in\mathcal{U}}\lambda\int_0^{+\infty} e^{-\lambda s}g(y(s,u,y_0),u(s))ds
 \end{equation}

In this framework the problem was initially to know whether, for a
given finite horizon $T$ and a given starting point $y_0$, a
minimising control $u$ existed, solution of the optimal control
problem $(T,y_0)$. Systems with large, but fixed horizons were
considered and, in particular, the class of "ergodic" systems (that
is, those in which any starting point in the state space $\Omega$ is
controllable to any point in $\Omega$) has been thoroughly studied.
These systems are asymptotically independent of the starting point
as the horizon goes to infinite. When the horizon is infinite, the
literature on optimal control has mainly focussed on properties of
given trajectories as the time tends to infinity. This approach
corresponds to the uniform approach in a game theoretical framework
and is often opposed to the asymptotic approach (described below),
which we have considered in what follows, and which has received
considerably less attention.

In a game-theoretical, discrete time framework, the same kind of
problem was considered since\cite{Sh}, but with several differences
in the approach: 1) the starting point may be chosen at random, i.e,
a probability $\mu$ may be given on $\Omega$, which randomly
determines the point from which the controller will start the play;
2) the controllability-ergodicity condition is generally not
assumed; 3) because of the inherent recursive structure of process
played in discrete time, the problem is generally considered for all
initial states and time horizons.

For these reasons, what is called the "asymptotic approach" --- the
behavior of $V_t(\cdot)$ as the horizon $t$ tends to infinity, or of
$V_{\lambda}(\cdot)$ as the discount factor $\lambda$ tends to zero
--- as been more studied in this discrete-time setup. Moreover, when it is considered in Optimal Control, in most cases \cite{Arisawa, Bettiol} an
ergodic assumption is made which not only ensures the convergence of
$V_t(y_0)$ to some $V$, but also forces the limit function $V$ to be
independent of the starting point $y_0$. The general asymptotic
case, in which no ergodicity condition is assumed, has been to our
knowledge studied for the first time recently. In \cite{Carda,QR}
the authors prove in different frameworks the convergence of
$V_t(\cdot)$ and $V_\lambda(\cdot)$ to some non-constant function
$V(y_0)$.

Some important, closely related questions are the following : does
the convergence of $V_t(\cdot)$ imply the convergence of
$V_{\lambda}(\cdot)$ ? Or vice versa ? If they both converge, does
the limit coincide ? A partial answer to these questions goes back
to the beginning of the $20^{th}$ century, when Hardy and Littlewood
proved (see \cite{HaLi}) that for any sequence of bounded real
numbers, the convergence of the Cesaro means is equivalent to the
convergence of their Abel means, and that the limits are then the
same :
%
\begin{theo}[Hardy-Littlewood 1914]
For any bounded sequence of reals $\{a_n\}_{n \geq 1}$, define $V_n=\frac{1}{n}\sum_{i=1}^n a_i$ and $V_\lambda=\lambda \sum_{i=1}^{+\infty} (1-\lambda)^{i-1} a_i$. Then,
\[ \underset{n\rightarrow +\infty}{\liminf}\ V_n\leq \underset{\lambda\rightarrow 0}{\liminf}\ V_\lambda\leq\underset{\lambda\rightarrow 0}{\limsup}\ V_\lambda\leq\underset{n\rightarrow +\infty}{\limsup}\ V_n. \]
Moreover, if the central inequality is an equality, then all inequalities are equalities.
\end{theo}
Noticing that $\{a_n\}$ can be viewed as a sequence of costs for
some deterministic (uncontrolled) dynamic in discrete-time, this
results gives the equivalence between the convergence of $V_t$ and
the convergence of $V_{\lambda}$, to the same limit. In 1971,
setting $V_t=\frac{1}{t}\int_{0}^t g(s)ds$ and $V_\lambda=\lambda
\int_{0}^{+\infty} e^{-\lambda s}g(s) ds$, for a given
Lebesgue-measurable, bounded, real function $g$, Feller proved in
\cite{Feller} that the same result holds for continuous-time
uncontrolled dynamics.
\begin{theo}[Feller 1971]
\[\underset{n\rightarrow +\infty}{\liminf}\ V_n\leq \underset{\lambda\rightarrow 0}{\liminf}\ V_\lambda\leq\underset{\lambda\rightarrow 0}{\limsup}\ V_\lambda\leq\underset{n\rightarrow +\infty}{\limsup}\ V_n.\]
Moreover, if the central inequality is an equality, then all inequalities are equalities.
\end{theo}
In 1992, Lehrer and Sorin \cite{LehSys} considered a discrete-time controlled dynamic, defined by a correspondence $\Gamma:\Omega\rightrightarrows\Omega$, with nonempty values, and by $g$, a bounded real cost function defined on $\Omega$. A feasible play at $z \in\Omega$ is an infinite sequence $\textbf{y}=\{y_n\}_{n\geq 1}$ such that $y_1=z$ and $y_{n+1}\in\Gamma(y_n)$. 
The value functions are defined by $V_n(z)=\inf \ \frac{1}{n}\sum_{i=1}^n g(y_i)$ and respectively $V_\lambda(y_0)=\inf \ \lambda \sum_{i=1}^{+\infty} (1-\lambda)^{i-1} g(y_i)$, where the infima are taken over the feasible plays at $z$.

\begin{theo}[Lehrer-Sorin 1992]
\[\underset{n\rightarrow +\infty}{\lim}\ V_n(z)=V(z) \text{ uniformly on } \Omega \Longleftrightarrow  \underset{\lambda\rightarrow 0}{\lim}\ V_\lambda(z)=V(z) \text{ uniformly on } \Omega.\]
\end{theo}
This result establishes the equivalence between uniform convergence
of $V_{\lambda}(y_0)$ when $\lambda$ tends to $0$ and uniform
convergence of $V_n(y_0)$ as $n$ tends to infinity, in the general
case where the limit may depend on the starting point $y_0$. The
uniform condition is necessary: in the same article, the authors
provide an example where only pointwise convergence holds and the
limits differs.

In 1998, Arisawa (see \cite{Arisawa}) considered a
continuous-time controlled dynamic and proved the equivalence
between the uniform convergence of $V_{\lambda}$ and the uniform
convergence of $V_t$ in the specific case of limits independent of
the starting point.

\begin{theo}[Arisawa 1998]\label{Arsw} Let $d\in\mathbb{R}$, then
\[\underset{t\rightarrow +\infty}{\lim}\ V_t(z)=d \text{ uniformly on } \Omega \Longleftrightarrow  \underset{\lambda\rightarrow 0+}{\lim}\ V_\lambda(z)=d \text{ uniformly on } \Omega.\]
\end{theo}
This does not settle the general case, in which the limit function
may depend on the starting point\footnote{Lemma 6 and Theorem 8 in
\cite{Arisawa} deal with this general setting, but we believe them
to be incorrect since they are stated for pointwise convergence and,
consequently, are contradicted by the example in Section
\ref{sectioncontreexemple}.}. For a continuous-time controlled
dynamic in which $V_t(y_0)$ converges to some function $V(y_0)$,
dependent on the state variable $y_0$, as $t$ goes to infinity, we
prove the following


\begin{theo}\label{OB-V}
$V_t(y_0)$ converges to $V(y_0)$ uniformly on $\Omega$, if and only
if $V_{\lambda}(y_0)$ converges to $V(y_0)$ uniformly on $\Omega$.
\end{theo}

In fact, we will prove this result in a more general framework, as
described in section 2. Some basic lemmas which occur to be
important tools will also be proven on that section. Section 3 will
be devoted to the proof of our main result. Section 4 will conclude
by pointing out, via an example, the fact that uniform convergence
is a necessary requirement for the Theorem \ref{OB-V} to hold. A
very simple dynamic is described, in which the pointwise limits of
$V_t(\cdot)$ and $V_{\lambda}(\cdot)$ exist but differ. It should be
noted that our proofs (as well as the counterexample in Section
\ref{sectioncontreexemple}) are adaptations in this continuous-time
framework of ideas employed in a discrete-time setting in
\cite{LehSys}. In the appendix we also point out that an alternate
proof of our theorem is obtained using the main theorem in
\cite{LehSys} as well as a discrete/continuous equivalence argument.

For completeness, let us mention briefly this other approach,
mentioned above as the uniform approach, and which has also been
deeply studied, see for exemple \cite{Carlson, Colonius, Grune}). In
these models, the optimal average cost value, i.e the $V_t$, is not
taken over a finite period of time $[0,t]$, which is then studied
for $t$ growing to infinite, as in \cite{LehSys, Arisawa, QR, HaLi,
Feller} or in our framework. On the contrary, only infinite
trajectories are considered, among which the value $\overline{V_t}$
is defined as $\inf_{u \in \mathcal{U}} \sup_{\tau \geq
t}\frac{1}{\tau}\int_0^{\tau}g(y(s,u,y_0),u(s))ds$, or some other
closely related variation. The asymptotic behavior, as $t$ tends to
infinity, of the function $\overline{V_t}$ has also been studied. In
\cite{Grune}, both $\lambda$-discounted and average evaluations of
an infinite trajectory are considered and their limits are compared.
However, we stress out that the asymptotic behavior of those
quantities is in general\footnote{The reader may verify that this is
indeed not the case in the example of Section
\ref{sectioncontreexemple}.} not related to the asymptotic behavior
of $V_t$ and $V_\lambda$.

\section{Model}
\subsection{General framework}\label{generalframework}
We consider a deterministic dynamic programming problem in
continuous time, defined by a measurable set of states $\Omega$, a
subset $\mathcal{T}$ of Borel-measurable functions from
$\mathbb{R_+}$ to $\Omega$, and a bounded Borel-measurable
real-valued function $g$ defined on $\Omega$. Without loss of
generality we assume $g:\Omega \to [0,1]$. For a given state $x$,
define $\Gamma(x):=\{X\in\mathcal{T},\ X(0)=x\}$ the set of all
feasible trajectories starting from $x$. We assume $\Gamma(x)$ to be
non empty, for all $x \in \Omega$.
Furthermore, the correspondence $\Gamma$ is closed under
concatenation: given a trajectory $X \in \Gamma(x)$ with $X(s)=y$,
and a trajectory $Y \in \Gamma(y)$, the concatenation of $X$ and $Y$
at time $s$ is
\begin{eqnarray}\label{concatenation}
 X \circ_s Y := \left\{
 \begin{array}{rl}
   X(t) & \text{if } t \leq s\\
   Y(t-s) & \text{if } t \geq s\\
 \end{array} \right.
\end{eqnarray}
\noindent and we assume that $X \circ_s Y \in  \Gamma(x)$.


We are interested in the asymptotic behavior of the average and the
discounted values. It is useful to denote the average payoff of a
play (or trajectory) $X \in \Gamma(x)$ by:
\begin{eqnarray}
\gamma_t(X)&:=&\dfrac{1}{t}\int_0^t g(X(s))ds\\
\gamma_\lambda(X)&:=&\lambda \int_0^{+\infty} e^{-\lambda
s}g(X(s))ds
\end{eqnarray}
This is defined for $t, \lambda \in (0, \infty)$. Naturally, we
define the values as:
\begin{eqnarray}
\label{eqVt}V_t(x)&=&\inf_{X\in\Gamma(x)} \gamma_t(X)\\
\label{eqVlambda}V_{\lambda}(x)&=&\inf_{X\in\Gamma(x)}
\gamma_{\lambda}(X)
\end{eqnarray}

Our main contribution is the following :

\begin{theo}\label{mainprop}
(A)  $V_\lambda \underset{\lambda \to 0}{\longrightarrow} V$,
uniformly  on  $\Omega$ $\Longleftrightarrow$ (B)  $V_t \underset{t
\to \infty}{\longrightarrow} V$, uniformly on
 $\Omega$.
\end{theo}

Notice that our model is a natural adaptation to the continuous-time
framework of deterministic dynamic programming problems played in
discrete time ; as it was pointed out during the introduction, this
theorem is an extension to the continuous-time framework of the main
result of \cite{LehSys}, and our proof use similar technics.

This result can be applied to the model presented in
section \ref{Controleoptimal}: denote
$\widetilde{\Omega}=\mathbb{R}^d\times U$ and for any
$(y_0,u_0)\in\widetilde{\Omega}$, define
$\widetilde{\Gamma}(y_0,u_0)=\{(y(\cdot),u(\cdot)) \ |
u\in\mathcal{U}, u(0)=u_0\text{ and } y\text{ is the solution of
(\ref{eqcontrole}).} \}$ Then $\widetilde{\Omega}$,
$\widetilde{\Gamma}$ and $g$ satisfy the assumptions of this
section. Defining $\widetilde{V}_t$ and $\widetilde{V}_\lambda$ as
in (\ref{eqVt}) and (\ref{eqVlambda}) respectively, since the
solution of (\ref{eqcontrole}) does not depend on $u(0)$ we get that
\begin{eqnarray*}
\widetilde{V}_t(y_0,u_0)&=&V_t(y_0)\\
\widetilde{V}_\lambda(y_0,u_0)&=&V_\lambda(y_0).
\end{eqnarray*}

Theorem \ref{mainprop} applied to $\widetilde{V}$ thus implies that
$V_t$ converges uniformly to a function $V$ in $\Omega$ if and only
if $V_\lambda$ converges uniformly to $V$ in $\Omega$.




\subsection{Preliminary results} We follow the ideas of
\cite{LehSys}, and start by proving two simple lemmas yet important
tools, that will be used in the proof. The first establishes that
the value increases along the trajectories. Then, we prove a
convexity result linking the finite horizon average payoffs and the
discounted evaluations on any given trajectory.

\begin{lemme}\label{decroissancetraj} Monotonicity (compare with Proposition 1 in \cite{LehSys}) \\
For all $X \in \mathcal{T}$, for all $s \geq 0$, we have
\begin{eqnarray}
 \liminf_{t \to \infty} V_t(X(0)) & \leq & \liminf_{t \to \infty} V_t(X(s)) \\
 \liminf_{\lambda \to 0} V_\lambda(X(0)) &\leq& \liminf_{\lambda \to 0} V_\lambda(X(s))
\end{eqnarray}
\end{lemme}
\begin{proof}
Set $y:=X(s)$ and $x:=X(0)$. For $\varepsilon>0$, take
$T\in\mathbb{R}_+$ such that $\frac{s}{s+T}<\varepsilon$. Let $t>T$
and take a $\varepsilon$-optimal trajectory for $V_t$, i.e. $Y\in
\Gamma(y)$ such that $\gamma_t(Y) \leq V_t(y)+\varepsilon$. Define
the concatenation of $X$ and $Y$ at time $s$ as in
(\ref{concatenation}), where $X \circ_s Y$ is in $\in\Gamma(x)$ by
our hypothesis. Hence
\begin{eqnarray*}
V_{t+s}(x) \leq \gamma_{t+s}(X \circ_s Y)&=& \dfrac{s}{t+s}
 \gamma_s(X) + \frac{t}{t+s} \gamma_{t}(Y) \\ &\leq& \varepsilon + \gamma_t(Y)\\
 &\leq& 2\varepsilon+V_t(y).
\end{eqnarray*}
Since this is true for any $t\geq T$ the result follows.

Similarly, for the discounted case let $\lambda_0>0$ be such that $
\lambda_0 \int_0^s e^{-\lambda_0 r} dr = 1-e^{\lambda_0 s}
<\varepsilon$. Let $\lambda \in (0,\lambda_0)$ and take $Y \in
\Gamma(y)$ a $\varepsilon$-optimal trajectory for $V_\lambda(y)$.
Then:
\begin{eqnarray*}
 V_\lambda(x)\leq\gamma_{\lambda}( X \circ_s Y )&=&  \lambda \int_0^{s}e^{-\lambda r}g(X(r))dr+ \lambda \int_s^{\infty}e^{-\lambda r}g(Y(r-s))dr\\
&\leq& \varepsilon + e^{-\lambda s}\gamma_{\lambda}(Y)\\
&\leq& 2\varepsilon + V_\lambda(y).
\end{eqnarray*}
Again, this is true for any $\lambda \in (0, \lambda_0)$, and the
result follows.

\end{proof}


\begin{lemme}\label{convexitylemma} Convexity (compare with equation (1) in \cite{LehSys})\\
For any play $X \in \mathcal{T}$, for any $\lambda>0$:
\begin{equation}\label{convexity}
 \gamma_{\lambda}(X)= \int_0^{\infty}\gamma_s(X)\mu_{\lambda}(s)ds,
\end{equation}
where $\mu_{\lambda}(s)ds:=\lambda^2 se^{-\lambda s }ds$ is a
probability density on $[0, +\infty]$.
\end{lemme}
\begin{proof}
It is enough to notice that the following relation holds, by
integration by parts :
\[ \gamma_{\lambda}(X)= \lambda \int_0^{\infty}e^{-\lambda s}g(X(s))ds = \lambda^2 \int_0^{\infty}se^{-\lambda s}\left(\frac{1}{s}\int_0^s g(X(r))dr \right)ds, \]
and that $\int_0^{\infty}\lambda^2 se^{-\lambda s }ds=1$.
\end{proof}

The probability measure $\mu_{\lambda}$ will play an important role
in the rest of the paper. Denoting $M(\alpha, \beta; \lambda):=
\displaystyle{\int_{\alpha}^{\beta}\mu_{\lambda}(s)ds}=e^{-\lambda
\alpha}(1+\lambda \alpha)-e^{-\lambda \beta}(1+\lambda \beta)$, we
prove here two estimates that will be helpful in the next section.
\begin{lemme}\label{poids}  
The two following results hold (compare with Lemma 3 in
\cite{LehSys}):
\begin{enumerate}
 \item[$\mathrm{(i)}$]$  \forall t>0, \exists \epsilon_0 \ such \ that \ \forall \epsilon \leq \epsilon_0,$ $M((1-\epsilon)t, t; 1/t)\geq \epsilon/{2e}.$

\item[$\mathrm{(ii)}$]$ \forall \delta>0, \exists \epsilon_0 \ such \ that \ \forall \epsilon \leq \epsilon_0$, $\forall t>0$, $M(\epsilon t, (1-\epsilon)t; 1/{t \sqrt{\epsilon}})\geq 1-\delta.$
\end{enumerate}
    \end{lemme}
\begin{proof}
Notice that in these particular cases, $M$ does not depend on $t$:
\begin{enumerate}
 \item[$\mathrm{(i)}$] $M(t(1-\epsilon),t;1/t)=(2-\epsilon)e^{-1+\epsilon}-2e^{-1}=\frac{1}{e}(\epsilon+o(\epsilon))\geq\frac{\epsilon}{2e}$, for $\epsilon$ small enough.

\item[$\mathrm{(ii)}$] $M(t\epsilon,t(1-\epsilon);1/\sqrt{\epsilon}t)=(1+\sqrt{\epsilon})e^{-\sqrt{\epsilon}}-(1-1/\sqrt{\epsilon}+\sqrt{\epsilon})e^{-1/\sqrt{\epsilon}+\sqrt{\epsilon}}$. This expression tends to $1$ as $\epsilon \to 0$, hence the result.
\end{enumerate}
\end{proof}

\section{Proof of Theorem \ref{mainprop}}

\subsection{From $V_t$ to $V_\lambda$}
Assume (B) : $V_t(\omega)$ converges to some $V(\omega)$ as $t$ goes
to infinity, uniformly on $\Omega$.

\begin{prop}\label{VlambdaVt+epsilon} (Compare with Proposition 4 in
\cite{LehSys})\\
 For all $\epsilon>0$, there exists $\lambda_0>0$
such that $V_{\lambda}(x)\geq V(x)-\varepsilon$ for every
$x\in\Omega$ and for all $\lambda \in (0,\lambda_0]$.
\end{prop}

\begin{proof}
Let $T$ be such that $\|V_t-V\|_\infty\leq\epsilon/2$ for every
$t\geq T$. Choose $\lambda_0>0$ such that $\displaystyle{ \lambda^2
\int_T^{\infty} {se^{-\lambda s}ds}=1-(1+\lambda T)e^{-\lambda
T}\geq 1-\epsilon/4}$, for every $\lambda \in (0,\lambda_0)$. Fix
$\lambda \in (0,\lambda_0)$ and take a play $Y  \in \Gamma(x)$ which
is $\epsilon/4$-optimal play  for $V_{\lambda}(x)$. Since
$\gamma_s(X) \geq 0$, the convexity formula (\ref{convexity}) from
Lemma \ref{convexitylemma} gives:
\begin{eqnarray*}
V_{\lambda}(x) +\epsilon/4 \geq \gamma_{\lambda}(Y) &\geq& 0 + \lambda^2 \int_T^{\infty}se^{-\lambda s }\gamma_s(Y)\\
&\geq& \lambda^2 \int_T^{\infty}se^{-\lambda s }V_s(x)\\
&\geq& (1-\epsilon/4)(V(x)-\epsilon/2)\\
&=& V(x)-\epsilon/4 V(x)-\epsilon/2+\epsilon^2/8\\
&\geq& V(x)-3/4\epsilon
\end{eqnarray*}
\end{proof}

\begin{lemme}\label{Lemme8LeSo}(Compare with Lemma 8 in
\cite{LehSys})\\
 $\forall \epsilon>0, \exists M$ such that for all $t \geq M, \forall x\in\Omega$, there is a play $X \in \Gamma(x)$ such that $\gamma_s(X) \leq V(x)+\epsilon$ for all $s \in [\epsilon t, (1-\epsilon)t].$
\end{lemme}
\begin{proof}
By (B) there exists $M$ such that $\|V_r-V\|\leq
\frac{\epsilon^2}{3}$ for all $r\geq \epsilon M$.
Given $t\geq M$ and $x\in\Omega$, let $X \in \Gamma(x)$ be a play
(from $x$) such that $\gamma_t(X)\leq V_t(x)+\frac{\epsilon^2}{3}$.
For any $s\leq(1-\epsilon)t$, $t-s\geq\epsilon t\geq \epsilon M$ so
Proposition \ref{decroissancetraj} (Monotonicity) imply that
\begin{equation}\label{eqVtmoinss}
V_{t-s}(X(s))\geq V(X(s))-\frac{\epsilon^2}{3}\geq
V(x)-\frac{\epsilon^2}{3}.
\end{equation}
Since $V(x)+\epsilon^2/3 \geq V_t(x)$, we also have:
\begin{eqnarray*}
t\left(V(x)+2\frac{\epsilon^2}{3}\right)&\geq&t\left(V_t(x)+\frac{\epsilon^2}{3}\right)\\
&\geq&t\gamma_t(X) = \int_0^s g(X(r))dr + \int_s^t g(X(r))dr \\
&\geq&s\gamma_s(X)+(t-s)V_{t-s}(X(s))\\
&\geq&s\gamma_s(X)+(t-s) \left(V(x)-\frac{\epsilon^2}{3}\right)
\text{ by (\ref{eqVtmoinss})}.
\end{eqnarray*}

Isolating $\gamma_s(X)$ we get:
\begin{eqnarray*}
 \gamma_s(X)&\leq&V(x)+\epsilon^2\frac{t}{s}\\
&\leq& V(x)+\epsilon, \quad \mathrm{for}\ s/\epsilon \geq t,
\end{eqnarray*}
and we have proved the result for all $s \in [\epsilon t, (1-
\epsilon)t]$.
\end{proof}

\begin{prop}\label{VlambdaV-delta}(Compare with Lemma 9 in
\cite{LehSys})\\
 $\forall \delta>0, \exists \lambda_0$ such that $\forall x\in\Omega$, for all $\lambda \in (0,\lambda_0]$, we have $V_{\lambda}(x)\leq V(x)+\delta$.
\end{prop}
\begin{proof}
By Lemma \ref{poids} (ii), one can choose $\epsilon$ small enough
such that $M(\epsilon t, (1-\epsilon)t; \frac{1}{t
\sqrt{\epsilon}})\geq 1-\frac{\delta}{2}$, for any $t$. In
particular, we can take $\epsilon \leq\frac{\delta}{2}$. Using Lemma
\ref{Lemme8LeSo} with $\delta/2$, we get that for $t \geq t_0$ (and
thus for $\lambda_t:=\frac{1}{t \sqrt{\epsilon}}\leq\frac{1}{t_0
\sqrt{\epsilon}}$) and for any $x\in\Omega$, there exists a play $X
\in \Gamma(x)$ such that
 \begin{eqnarray*}
 \gamma_{\lambda_t}(X)&\leq& \delta/2 + {\lambda_t}^2\int_{\epsilon t}^{(1-\epsilon)t} s e^{\lambda_t s} \gamma_s(X) ds\\
 &\leq& \delta/2 + 1 \cdot (V(x)+\delta/2).
 \end{eqnarray*}
\end{proof}

Propositions \ref{VlambdaVt+epsilon} and \ref{VlambdaV-delta}
establish the first part of Theorem \ref{mainprop}: $(B) \Rightarrow
(A)$.

\subsection{From $V_\lambda$ to $V_t$}
Now assume $(A)$ : $V_{\lambda}(\omega)$ converges to some
$W(\omega)$ as $\lambda$ goes to 0, uniformly on $\Omega$. We start
by a technical Lemma:

\begin{lemme}\label{localaverages}(Compare with Proposition 2 in
\cite{LehSys})\\
Let $\epsilon>0$. For all $x\in \Omega$ and $t > 0$, and for any
trajectory $Y \in \Gamma(x)$ which is $\epsilon/2$-optimal for the
problem with horizon $t$, there is a time $L \in
[0,t(1-\epsilon/2)]$ such that, for all $T \in ]0,t-L]$:
\[ \frac{1}{T}\int_L^{L+T} g(Y(s))ds \leq V_t(x)+\epsilon. \]
\end{lemme}
\begin{proof}
Fix $Y\in \Gamma(x)$ some $\epsilon/2$-optimal play for $V_t(x)$.
The function $s\rightarrow\gamma_s(Y)$ is continuous on $]0,t]$ and
satisfies $\gamma_t(Y) \leq V_t(x)+\epsilon/2$. The bound on $g$
implies that $\gamma_r(Y) \leq V_t(x)+\epsilon$ for all $r \in
[t(1-\epsilon/2),t]$.

Consider now the set $\{s\in ]0,t] \ | \ \gamma_s(Y)>V_t(x)+\epsilon
\}$.
If this set is empty, then take $L=0$ and observe that
$\frac{1}{r}\int_0^{r} g(Y(s))ds \leq V_t(x)+\epsilon, \forall r \in
]0,t]$.

Otherwise, let $L$ be the superior bound of this set. Notice that
$L<t(1-\epsilon/2)$ and that by continuity
$\gamma_L(Y)=V_t(x)+\epsilon$. Now, for any $T \in [0,t-L]$,
\begin{eqnarray*}
V_t(x)+\epsilon&\geq& \gamma_{L+T}(Y)\\
&=&\frac{L}{L+T}\gamma_L(Y) +
\frac{T}{L+T}\left(\frac{1}{T}{\int_L^{L+T} g(Y(s))ds}\right)\\
&=&\frac{L}{L+T}\left(V_t(x)+\epsilon
\right)+\frac{T}{L+T}\left(\frac{1}{T}{\int_L^{L+T}
g(Y(s))ds}\right)
\end{eqnarray*}
and the result follows.
\end{proof}

\begin{prop}\label{VtWepsilon}(Compare with Lemma 6 in
\cite{LehSys})\\
 $\forall \epsilon>0, \exists T$ such that for all $t \geq T$ we have $V_t(x) \geq W(x)-\epsilon$, for all $x \in \Omega.$
\end{prop}
\begin{proof} Let $\lambda$ be such that $\|V_{\lambda}-W\| \leq \epsilon/8$, and $T$ such that $\displaystyle{ \lambda^2 \int_{T \epsilon/4}^{\infty}se^{-\lambda s }ds} <\epsilon/8$. Proceed by contradiction and suppose that $\epsilon>0$ is such that for every $T$, there exists $t_0 \geq T$ and a state $x_0 \in \Omega$ such that $V_{t_0}(x_0)<W(x_0)-\epsilon$.

Using Lemma \ref{localaverages} with $\epsilon/2$, we get a play $Y
\in \Gamma(x_0)$ and a time $L \in [0,t_0(1-\epsilon/4)]$ such that,
$\forall s \in [0, t_0-L]$ (and, in particular, $\forall s \in [0,
t_0 \epsilon /4]$):
\[\frac{1}{s}\int_L^{L+s} g\left(Y(r)\right)dr \leq V_{t_0}(x_0)+ \epsilon/2 < W(x_0)-\epsilon/2.\]
Thus,
\begin{eqnarray*}
W(Y(L))-\epsilon/8  &\leq&  V_\lambda(Y(L))\\
& \leq & \lambda \int_0^{\infty}e^{-\lambda s}g(Y(L+s))ds\\
&\leq & \lambda^2 \int_0^{t_0 \epsilon/4}se^{-\lambda s }\left(\frac{1}{s}\int_L^{L+s} g\left(Y(r)\right)dr\right)ds + \ \epsilon/8\\
&\leq& W(x_0)-\epsilon/2+\epsilon/8\\
& = & W(x_0)-3\epsilon/8.
\end{eqnarray*}
This gives us $W(Y(L))\leq W(x_0)-\epsilon/4$, contradicting
Proposition \ref{decroissancetraj} (Monotonicity).
\end{proof}

\begin{prop}\label{VtW-epsilon}(Compare with Lemma 7 in
\cite{LehSys})\\
 $\forall \epsilon>0, \exists T$ such that for all $t \geq T$ we have $V_t(x) \leq W(x)+\epsilon$, for all $x \in \Omega$.
\end{prop}

\begin{proof}
 Otherwise, $\exists \epsilon>0$ such that $\forall T, \ \exists t \geq T$ and $x\in \Omega$ with $V_t(x)>W(x)+\epsilon.$ For any $X \in \Gamma(x)$ consider the (continuous in $s$) payoff function $\gamma_s(X)=\displaystyle{\frac{1}{s}\int_0^s g(X(r))dr}$. Of course, $\gamma_t(X)\geq V_t(x)>W(x)+\epsilon$. Furthermore, because of the bound on $g$,
\[\gamma_r(X)\geq W(x)+\epsilon/2, \ \forall r\in [t(1-\epsilon/2), t].\]


By Lemma \ref{poids}, we can take $\epsilon$ small enough, so that
for all $t$, $M(t(1-\epsilon/2),t;1/t)\geq \epsilon/4e$ holds. We
set $\delta:=\epsilon/4e$.

By Proposition \ref{VtWepsilon}, there is a $K$ such that $V_t \geq
W(x)-\delta \epsilon/8$, for all $t \geq K$.

For $K$ fixed, we consider $M(0,K;1/t)=1-e^{-K/t}(1+K/t)$ as a
function of $t$. Clearly, it tends to $0$ as $t$ tends to infinity,
so let $t$ be such that this quantity is smaller than $\delta
\epsilon/16$. Also, let $t$ be big enough so that $\Arrowvert
V_{1/t}-W\Arrowvert<\delta \epsilon/5$, which is a consequence of
assumption (A).

On the following, we set $\tilde{\lambda}:=1/t$ and consider the
$\tilde{\lambda}$-payoff of some play $X \in \Gamma(x)$. We'll split
the integral over $[0,+\infty]$ in three parts : $\mathcal{K}=[0,K],
\mathcal{R}=[t(1-\epsilon/2),t]$, and $(\mathcal{K} \cup
\mathcal{R})^c$. The three parts are clearly disjoint since $t>>K$.
We have seen that $\mu_{\lambda}(s)ds=\lambda^2se^{-\lambda s}ds$ is
a probability measure on $[0,+\infty]$, for any $\lambda>0$. Then by
the Convexity formula (\ref{convexity}), we can write:
\[\gamma_{\tilde{\lambda}}(X)=\left(\int_{\mathcal{K}}\gamma_s(X)\mu_{\tilde{\lambda}}(ds)+\int_{\mathcal{R}}\gamma_s(X)\mu_{\tilde{\lambda}}(ds)+\int_{(\mathcal{K} \cup \mathcal{R})^c}\gamma_s(X)\mu_{\tilde{\lambda}}(ds) \right).\]

Recall that \begin{eqnarray*}
\gamma_s(X)_{|\mathcal{K}} &\geq& 0\\
\gamma_s(X)_{|(\mathcal{K} \cup \mathcal{R})^c} &\geq& W(x)-\delta \epsilon/8\\
\gamma_s(X)_{|\mathcal{R}} &\geq& W(x)+\epsilon/2
            \end{eqnarray*}

It is straightforward that
\begin{eqnarray*}
\gamma_{\tilde{\lambda}}(X) &\geq& 0 + \delta \cdot (W(x)+\epsilon/2) +(1-\delta-\delta \epsilon/16) \cdot (W(x)-\delta \epsilon/8)\\
&=&W(x)+ \delta \epsilon \left(\frac{1}{2}-\frac{1}{16}-\frac{1}{8}-\frac{\delta}{8}+\frac{\delta \epsilon}{64}\right)\\
& \geq &W(x)+\delta \epsilon /4.
\end{eqnarray*}
This is true for any play, so its infimum also satisfies
$V_{\tilde{\lambda}}(x)\geq W(x)+\delta \epsilon /4$, which is a
contradiction, for we assumed that $V_{\tilde{\lambda}}<W(x)+\delta
\epsilon /5$.
\end{proof}

Propositions \ref{VtWepsilon} and \ref{VtW-epsilon} establish the
second half of Theorem \ref{mainprop} : $(A) \Rightarrow (B)$.



\section{A counter example for pointwise
convergence}\label{sectioncontreexemple} In this section we give an
example of an optimal control problem in which both $V_t(\cdot)$ and
$V_\lambda(\cdot)$ converge pointwise on the state space, but to two
different limits. As implied by Theorem \ref{mainprop}, the
convergence is not uniform on the state space.

Lehrer and Sorin were the first to construct such an
example\cite{LehSys}, in the discrete-time framework. We consider
here one of its adaptation in continuous time, which was studied as
Example 5 in \cite{QR}\footnote{We thank Marc Quincampoix for
pointing out this example to us, which is simpler that our original
one.}, where the notations are the same that in Section
\ref{Controleoptimal}:
\begin{itemize}
\item The state space is $\Omega=\mathbb{R}_+^2$.
\item The payoff function is given by $g(x,y)=0$ if $x\in[1,2]$, 1 otherwise.
\item The set of control is $U=[0,1]$.
\item The dynamic is given by $f(x,y,u)=(y,u)$ (thus $\Omega$ is forward
invariant.)
\end{itemize}
An interpretation is that the couple $(x(t),y(t))$ represents the position and the speed of some mobile moving along an axis, and whose acceleration $u(t)$ is controlled. Observe that since $U=[0,1]$, the speed $y(t)$ increases during a play.
We claim that for any $(x_0,y_0)\in\mathbb{R}_+^2$, $V_t(x_0,y_0)$ (resp $V_\lambda(x_0,y_0)$) converges to $V(x_0,y_0)$ as $t$ goes to infinity (resp. converges to $W(x_0,y_0)$ as $\lambda$ tends to 0, where:
\begin{eqnarray*}
V(x_0,y_0)&=&\begin{cases}
1& \text{if } y_0>0 \text{ or } x_0>2\\
0& \text{if } y_0=0 \text{ and } 1\leq x_0\leq 2\\
\frac{1-x_0}{2-x_0} &\text{if } y_0=0 \text{ and } x_0< 1
\end{cases}\\
W(x_0,y_0)&=&\begin{cases}
1& \text{if } y_0>0 \text{ or } x_0>2\\
0& \text{if } y_0=0 \text{ and } 1\leq x_0\leq 2\\
1-\frac{(1-x_0)^{1-x_0}}{(2-x_0)^{2-x_0}} &\text{if } y_0=0 \text{ and } x_0<1
\end{cases}
\end{eqnarray*}
Here we only prove that $V(0,0)=\frac{1}{2}$ and
$W(0,0)=\frac{3}{4}$ ; the proof for $y_0=0$ and $0<x_0<1$ is
similar and the other cases are easy.

First of all we prove that for any $t$ or $\lambda$ and any
admissible trajectory (that is, any function $X(t)=(x(t),y(t))$
compatible with a control $u(t)$), starting from $(0,0)$,
$\gamma_t(X)\geq \frac{1}{2}$ and $\gamma_\lambda(X)\geq
\frac{3}{4}$. This is clear if $x(t)$ is identically 0, so consider
this is not the case. Since the speed $y(t)$ is increasing, we can
define $t_1$ and $t_2$ as the time at which $x(t_1)=1$ and
$x(t_2)=2$ respectively, and moreover we have $t_2\leq 2t_1$. Then,
\begin{eqnarray*}
\gamma_t(X)&=&\frac{1}{t}\left(\int_0^{min(t,t_1)}ds+\int_{min(t,t_2)}^t
ds\right)\\
&=&1+\min\left(1,\frac{t_1}{t}\right)-\min\left(1,\frac{t_2}{t}\right)\\
&\geq&1+\min\left(1,\frac{t_2}{2t}\right)-\min\left(1,\frac{t_2}{t}\right)\\
&\geq&\frac{1}{2}
\end{eqnarray*}
\noindent and
\begin{eqnarray*}
\gamma_\lambda(X)&=&\int_0^{t_1}\lambda e^{-\lambda
s}ds+\int_{t_2}^{+\infty} \lambda e^{-\lambda s}ds\\
&=&1-e^{-\lambda t_1}+e^{-\lambda t_2}\\
&\geq&1-e^{-\lambda t_1}+e^{-2\lambda t_1}\\
&\geq&\min_{a>0}\{1-a+a^2\}\\
&\geq& \frac{3}{4}.
\end{eqnarray*}
On the other hand, one can prove\cite{QR} that $\limsup V_t(0,0)\leq
\frac{1}{2}$ : in the problem with horizon $t$, consider the control
"$u(s)=1$ until $s=\frac{2}{t}$ and then 0". Similarly one proves
that $\limsup V_\lambda(0,0)\leq \frac{3}{4}$: in the
$\lambda$-discounted problem, consider the control "$u(s)=1$ until
$s=\frac{\lambda}{\ln2}$ and then 0".

So the functions $V_t$ and $V_\lambda$ converge pointwise on
$\Omega$, but their limits $V$ and $W$ are different, since we have just shown $V(0,0)\neq W(0,0)$. One can verify that neither
convergence is uniform on $\Omega$ by considering
$V_t(1,\varepsilon)$ and $V_\lambda(1,\varepsilon)$ for small
positive $\varepsilon$.

\begin{remarque}
One may object that this example is not very regular since the
payoff $g$ is not continuous and the state space is not compact.
However a related, smoother example can easily be constructed:
\begin{enumerate}
\item The set of controls is still $[0,1]$.
\item The continuous cost $g(x)$ is equal to 1 outside the segment [0.9,2.1],
to 0 on [1,2], and linear on the two remainings intervals.
\item  The compact state space is $\Omega=\{(x,y)|0\leq y \leq \sqrt{2x}\leq 2\sqrt{2}\}$.
\item The dynamic is the same that in the original example for $x\in[0,3]$, and $f(x,y,u)=((4-x)y,(4-x)u)$ for $3\leq x \leq 4$.
The inequality $y(t)y'(t)\leq x'(t)$ is thus satisfied on any
trajectory, which implies that $\Omega$ is forward invariant under
this dynamic.
\end{enumerate}
With these changes the values $V_t(\cdot)$ and $V_\lambda(\cdot)$
still both converge pointwise on $\Omega$ to some
$\widetilde{V}(\cdot)$ and $\widetilde{W}(\cdot)$ respectively, and
$\widetilde{V}(0,0)\neq\widetilde{W}(0,0)$.
\end{remarque}

\section{Possible extensions}
\begin{itemize}
\item We considered the finite horizon problem and the discounted
one, but it should be possible to establish similar Tauberian
theorems for other, more complex, evaluations of the payoff. This
was settled in the discrete time case in \cite{MonSys}.
\item It would be very fruitful to establish necessary or
sufficient conditions for uniform convergence to hold. In this
direction we mention \cite{QR} in which sufficient conditions for
the stronger notion of Uniform Value (meaning that there are
controls that are nearly optimal no matter the horizon, provided it
is large enough) are given in a general setting.
\item In the discrete case an example is constructed in \cite{MonSys} in which there is no uniform
value despite uniform convergence of the families $V_t$ and
$V_\lambda$. It would be of interest to construct such an example in
continuous time, in particular in the framework of section
\ref{Controleoptimal}.
\item It would be very interesting to study Tauberian theorems for dynamic systems that are controlled by two conflicting controllers. In the framework of differential games this has been done recently (Theorem 2.1 in \cite{AlvBar}): an extension of Theorem \ref{Arsw} has been accomplished for two player games
in which the limit of $V_T$ or $V_\lambda$ is assumed to be
independent of the starting point. The similar result in the
discrete time framework is a consequence of Theorems 1.1 and 3.5 in
\cite{KoNe}. Existence of Tauberian theorems in the general setup of
two-persons zero-sum games with no ergodicity condition remains open
in both the discrete and the continuous settings.
\end{itemize}

\section*{Acknowledgments}
This article was done as part of the PhD of the first author. Both
authors wish to express their thanks to Sylvain Sorin for his
numerous comments and his great help. We also thank Helena
Frankowska and Marc Quincampoix for helpful remarks on earlier
drafts.

\section*{Appendix}
We give here another proof\footnote{We thank Fr\'ed\'eric Bonnans for
the idea of this proof} of Theorem \ref{mainprop} by using the
analoguous result in discrete time \cite{LehSys} as well as an
argument of equivalence between discrete and continuous dynamic.

%

Consider a deterministic dynamic programming problem in continuous
time as defined in section \ref{generalframework}, with a state
space $\Omega$, a payoff $g$ and a dynamic $\Gamma$. Recall that,
for any $\omega \in \Omega$, $\Gamma(\omega)$ is the non empty set
of feasible trajectories, starting from $\omega$. We construct an
associated deterministic dynamic programming problem in
\emph{discrete} time as follows.

Let $\widetilde{\Omega}=\Omega\times[0,1]$ be the new state space and let
$\widetilde{g}$ be the new cost function, given by
$\widetilde{g}(\omega,x)=x$. We define a multivalued-function with
nonempty values $\widetilde{\Gamma}:\widetilde{\Omega}\rightrightarrows\widetilde{\Omega}$
by
\[
(\omega,x)\in\widetilde{\Gamma}(\omega',x')\Longleftrightarrow
\exists X \in \Gamma(\omega'),\text{ with } X(1)=\omega\text{ and } \int_0^1
g(X(t))dt=x.
\]
Following \cite{LehSys}, we define, for any initial state
$\widetilde{\omega}=(\omega, x)$
\begin{eqnarray}
v_n(\widetilde{\omega})&=&\inf \frac{1}{n} \sum_{i=1}^{n}
\widetilde{g}(\widetilde{\omega}_i)\\
v_\lambda(\widetilde{\omega})&=&\inf  \lambda \sum_{i=1}^{+\infty}
(1-\lambda)^{i-1}\widetilde{g}(\widetilde{\omega}_i)
\end{eqnarray}
\noindent where the infima are taken over the set of sequences
$\{\widetilde{\omega}_i\}_{i\in\mathbb{N}}$ such that
$\widetilde{\omega}_0=\widetilde{\omega}$ and
$\widetilde{\omega}_{i+1}\in\widetilde{\Gamma}(\widetilde{\omega}_{i})$
for every $i \geq 0$.

Theorem \ref{mainprop} is then the consequence of the following
three facts. Firstly, the main theorem of Lehrer and Sorin in \cite{LehSys}, which states that uniform convergence (on $\widetilde{\Omega}$)
of $v_n$ to some $v$ is equivalent to uniform
convergence of $v_\lambda$ to the same $v$.

Secondly, the concatenation hypothesis (\ref{concatenation}) on $\Gamma$ implies that
for any $(\omega,x)\in\widetilde{\Omega}$
\[
v_n(\omega,x)=V_n(\omega)
\]
\noindent where $V_t(\omega)=\inf_{X\in\Gamma(\omega)} \frac{1}{t}
\int_0^n g(X(s)) ds$, as defined in equation
(\ref{eqVt}). Consequently, because of the bound on $g$, for any
$t\in\mathbb{R}+$ we have
\[
|V_t(\omega)-v_{\lfloor t \rfloor}(\omega,x)|\leq \frac{2}{\lfloor t
\rfloor}
\]
\noindent where ${\lfloor t \rfloor}$ stands for the integral part of $t$.

Finally, again because of hypothesis (\ref{concatenation}), for any
$\lambda\in ]0,1]$,
\[
v_\lambda(\omega,x)=\inf_{X\in\Gamma(\omega)}
\lambda\int_0^{+\infty} (1-\lambda)^{\lfloor t \rfloor} g(X(t)) dt.
\]
Hence, by equation (\ref{eqVlambda}) and the bound on the cost function, for
any $\lambda\in ]0,1]$,
\[
|V_\lambda(\omega)-v_\lambda(\omega,x)|\leq \lambda
\int_0^{+\infty} \left| (1-\lambda)^{\lfloor t \rfloor}-e^{-\lambda
t}\right|dt
\]
\noindent which tends uniformly (with respect to $x$ and $\omega$)
to 0 as $\lambda$ goes to 0 by virtue of the following lemma.
\begin{lemme}
The function \[\lambda\longrightarrow\lambda \int_0^{+\infty} \left|
(1-\lambda)^{{\lfloor t \rfloor}}-e^{-\lambda t}\right|dt\]
\noindent converges to 0 as $\lambda$ tends to 0
\end{lemme}

\begin{proof}
Since $\lambda \int_0^{+\infty} (1-\lambda)^{{\lfloor t\rfloor}} = \lambda \int_0^{+\infty}e^{-\lambda t}dt = 1$, for any $\lambda >0$, the lemma is equivalent to the
convergence to $0$ of
\[E(\lambda):=\lambda \int_0^{+\infty} \left[ (1-\lambda)^{{\lfloor t \rfloor}}-e^{-\lambda t}\right]_+dt\] \noindent where $[x]_+$ denotes the positive part of $x$.
Now, from the relation $1-\lambda \leq e^{-\lambda}$, true for any $\lambda$, one can easily deduce that, for any $\lambda >0, \ t\geq 0$, the relation $(1-\lambda)^{\lfloor t \rfloor}e^{\lambda t} \leq e^{\lambda}$ holds. Hence,

\begin{eqnarray*}
E(\lambda)&=&\lambda \int_0^{+\infty} e^{-\lambda t}\left[(1-\lambda)^{\lfloor t \rfloor}e^{\lambda t}-1 \right]_+dt\\
&\leq&\lambda \int_0^{+\infty} e^{-\lambda t}(e^{\lambda}-1) dt\\
&=&e^{\lambda}-1
\end{eqnarray*}
\noindent which converges to 0 as $\lambda$ tends to 0.
\end{proof}

 \end{document}